\documentclass[times,doublespace]{amsart}
\def\rev{\usepackage[active]{srcltx}}\rev
\usepackage[mathscr]{eucal}
\usepackage{amssymb, amsmath,array, amscd}
\usepackage{enumerate}
\usepackage[colorlinks]{hyperref}
\usepackage[all]{xy}

\def\p#1{\ensuremath{\mathbb{P}^{#1}}}

\def\ie{{\em i.e., }}
\def\C{\ensuremath{\mathbb C}}
\def\ve{\ensuremath{^\vee}}
\def\lar{\ensuremath{\,\longrightarrow\,}}
\def\be{\begin{equation}}\def\ee{\end{equation}}\def\bsm{\left(\begin{smallmatrix}}\def\esm{\end{smallmatrix}\right)}
\def\ba#1{\begin{array}{#1}}\def\ea{\end{array}}

\def\vez {\ensuremath{\times}}
\def\ra{\ensuremath{\,\rightarrow\,}}
\def\us{\ensuremath{^\star}}

\def\na#1{\noalign{\vskip#1pt}}

\newtheorem{coisa}{}[section]
\def\coi#1{\begin{coisa}{\bf
      #1.} \em}
\def\eco{\end{coisa}}
\def\cois#1{\begin{coisa}{\bf #1.}
\vspace{-10pt}
\eco
\addcontentsline{toc}{subsection}
{\numberline{}{#1}}
}

\def\runningheads#1#2{\markboth{\MakeTextUppercase{#1}}{\MakeTextUppercase{#2}}}

\def\rk{\ensuremath{\operatorname{rank}}}

\def\s#1{\ensuremath{\operatorname{Sym}_{#1}}}

\def\w{\ensuremath{\omega}}

\def\E{\ensuremath{\mathcal{V}}}

\newcommand{\Y}{\ensuremath{\mathbb Y}}

\newcommand{\grass}{{\mathbb G}}

\newcommand{\de}{\partial}
\newcommand{\OX}{\mathscr{O}}
\newcommand{\FF}{{\mathcal F}}

\newcommand{\T}{{\mathcal T}}
\newcommand{\Ker}{{\ensuremath{\mathrm{ker}}}}

\newcommand{\im}{{\ensuremath{\mathrm{Im}}}}

\def\mb#1{\ensuremath{\mathbf{#1}}}

\def\mn{\medskip\noindent}

\begin{document}

\message{ !name(pbcomponents.tex) !offset(-3) }


\runningheads{V. Ferrer and I. Vainsencher}{pullback components of the space of codimension one foliations}
\title[Linear Pullback components of the space of codimension one foliations]{linear pullback components of the space of codimension one foliations}
\author{V. Ferrer and I. Vainsencher}
\address{
Instituto de  Matem\'{a}tica e Estat\'istica, UFF -- Rua Prof. Marcos Waldemar de Freitas Reis, S/n, Bloco H, 5 andar. 24210-201, Niter\'{o}i, RJ - Brazil. 
 \newline\indent
ICEX--Depto. Matem\'{a}tica, UFMG -- Av. Antonio
Carlos 6627. 31270-901 Belo Horizonte Brazil. 
}

\thanks{The authors were partially supported by  CNPQ.}


\begin{abstract}
The space of holomorphic
foliations of codimension one and degree $d\geq 2$  in $\p n$ 
($n\geq 3$) has an irreducible
component whose general element can be written as a
pullback $F^*\FF$, where $\FF$ is a general 
foliation of degree $d$ in $\p 2$  and $F:\p
n\dashrightarrow \p 2$ is a general rational linear map. 
We give a polynomial formula  for the
degrees  of such components.

\end{abstract}

\maketitle

\section*{Introduction}\label{intro}
Codimension one holomorphic foliations in \p n 
are defined by  nonzero integrable twisted
$1-$forms, \ie   $\w\in H^0(\p n, \Omega_{\p n}(d+2) )$
satisfying
\begin{equation}\label{int}
\w\wedge d\w=0
\end{equation}

Since $\w,\lambda\w$ yield the same foliation for any
 $\lambda\in \C^*$, the space of such foliations  is in fact a closed
 subscheme  $\mathbb{F}(d,n)$  of the projective space
 $\p{} (H^0(\p n, \Omega_{\p n}(d+2) ))$ defined by the
 equations derived from  (\ref{int}). 

Explicitly, these equations are as follows.
Any element $\w$ of $H^0(\p n, \Omega_{\p n}(d+2)
)$
can be written
$\w=A_{0}dZ_{0}+\dots +A_{n}dZ_{n}$, where $A_{i}$ are homogeneous polynomials of degree  $d+1$ in the variables $Z_{0},\dots,Z_{n}$, such that 
$A_{0}Z_{0}+\dots +A_{n}Z_{n}=0$. 
The integrability condition (\ref{int}) imposes 
relations 
arising from
$$
A_{i}(\frac{\de A_{k}}{\de Z_{j}}-\frac{\de A_{j}}{\de Z_{k}})+A_{j}(\frac{\de A_{i}}{\de Z_{k}}-\frac{\de A_{k}}{\de Z_{i}})+A_{k}(\frac{\de A_{j}}{\de Z_{i}}-\frac{\de A_{i}}{\de Z_{j}})=0.
$$
 with  $0\leq i<j<k\leq n$.
 These equations are quadratic in the coefficients
 of the polynomials $A_{i}$. 
For  $n=2$
the space
 of foliations is a projective space:   the integrability condition is
 automatically satisfied.

The geometry of the space of codimension one
foliations in $\p n$ for $n\geq 3$ is a rich field
of research. In particular the problem of
describing the irreducible components of these
spaces has 
received many contributions,
cf.\, \cite{CLN}, \cite{CP}, \cite{CIP}, \cite{LPT} and
references therein just to quote a few.

In \cite{J} Jouanolou describes the irreducible
components for the space of foliations of degrees
$d=0,1$.  We have that $\mathbb{F}(0,n)$ is naturally
isomorphic to the grassmannian of subspaces of codimenion
2 in \p n.  The space of foliations of degree $d=1$ has
two irreducible components.

For foliations of degree $2$ in $\p n$, $n\geq 3$
Cerveau and Lins Neto \cite{CLN} have shown that
there are just $6$ components. We recall some
known components of $\mathbb{F}(d,n)$ pinpointing
those for which the degree has been found.
\begin{enumerate}
 
\item Logarithmic components. Let
$F_{1},\dots,F_{m}$ denote homogeneous polynomials
in $\C^{n+1}$, $m\geq 3$. Set $d_{i}=\deg(F_{i})$,
and let $\lambda_{1},\dots,\lambda_{m}\in \C^*$ be
such that $\sum_{i}\lambda_{i} d_{i}=0$. Then the
form $\w=F_{1}\dots
F_{m}\sum \lambda_{i}\frac{dF_{i}}{F_{i}}$ is an
integrable $1$-form in $\p n$. The foliation
induced by \w\ is called a logarithmic
foliation. Suppose that $F_{1},\dots,F_{m}$ are
irreducible and relatively prime. The space of
this type of foliations is denoted by
$\mathcal{L}(d_{1},\dots,d_{m})$.  Its closure
$\overline{\mathcal{L}(d_{1},\dots,d_{m})}$ is an
irreducible component of $\mathbb{F}(d,n)$, where
$d=d_{1}+\dots +d_{m}-2$, cf.\,\cite{CalA}. In the
case $d=2$ we have two possibilities:
$\overline{\mathcal{L}(1,1,1,1)}$ and
$\overline{\mathcal{L}(1,1,2)}$. Their degrees
have not been found to the best of our knowledge.

\medskip

\item Rational components. Let  
$F,G$ be homogeneous polynomials such that
\begin{itemize}

\item $\deg(F)=l , \deg(G)=m ,  m+l-2=d$;
\item $dF\wedge dG\neq 0$ along
$\{F=G=0\}\setminus \{0\}$ and
\item the hypersurface $\{F=0\}\subset \p n$ is smooth.
\end{itemize}
Then $\w=mGdF-lFdG$ defines a foliation in $\p n$
of degree $d$. Denote by $R(l,m)$ the set of
foliations of this type.  The closure
$\overline{R(l,m)}$ is an irreducible component of
$\mathbb{F}(d,n)$, cf.\,G\'omez-Mont and Lins Neto\,\cite{GMLN}.  Observe that
for $d=2$ we have two possibilities: $l=m=2$ or
$l=1,m=3$. In this way we obtain two more
components of $\mathbb{F}(2,n)$:
$\overline{R(2,2)}$ and $\overline{R(1,3)}$. The
degrees of these components were computed
in Cukierman-Pereira-Vainsencher\,\cite{CIP}.

\medskip

\item The exceptional component  $\mathcal{E}$ of
$\mathbb{F}(2,3)$. Here the leaves of a typical
foliation are the orbits of a linear action of the
affine group in one variable in $\p 3$.  It is
shown in \cite{CLN} that for $n>3$, in addition to
the components described above there exists a
component $\mathcal{E}_{n}$ obtained by linear
pullback $\p n \dashrightarrow \p 3$ of foliations
in $\mathcal{E}$. The degree of \
$\mathcal{E}$ \,has been found in
\cite{RV}.

\medskip

\item
Pullback components.  Let $\FF$ be a foliation of
degree $k$ in $\p 2$, defined by a $1$-form $\w$,
and $F: \p {n}\dashrightarrow \p 2$ a dominant
rational map of degree $m$. Then the $1$-form
$F^*\w$ defines a foliation $\FF^*$ of degree $d =
(k+2)m-2$ in $\p n$.  In \cite{CLNE} Cerveau,
Lins Neto and Edixhoven 
prove that for $n\geq 3$ the closure of the set of foliations
obtained by such  $F^*\w$ form an irreducible
component $PB(m,k,n)$ of $\mathbb{F}((k+2)m-2,n)$.

Our goal in this work is to
compute the degrees and dimensions of the
linear pullback components,
$LPB(d,n):=PB(1,d,n)$. Note that for
$d=(k+2)m-2\leq 6$ the only pullback components in
$\mathbb{F}(d,n)$ are linear or quadratic.  For
$d=7$ there exists a component of cubic pullbacks.
For $d=8, 9$ the only pullback components in
$\mathbb{F}(d,n)$ are linear or quadratic.  If
$d+2$ is prime, then the only pullback component
of $\mathbb{F}(d,n)$ is linear pullback.  Our main
result is a polynomial formula for $\deg
LPB(d,n)$, see \ref{formulas}.
\end{enumerate}

\section{Pull-back of projective 1-forms.}\label{pullbacksection}

\coi{Projective 1-forms}
The main reference for this material
is  \cite{J}.
A projective 1-form of
degree $d$ in $\p n$ is a global section of
$\Omega^1_{\p n}(d+2)$, for
some $d\geq 0$.  

We denote by $S_d$ the space 
$ 
H^0(\p n,\OX_{\p n}(d))={\rm Sym}_d(\C^{n+1})\ve
$ of
homogeneous polynomials of degree $d$ in the
variables $Z_0,\dots,Z_n$.  We write
$\partial_i=\partial/ \partial Z_i $  thought
of as a vector field basis for $\C^{n+1}$. The dual
basis will also be written as $dZ_0,\dots,dZ_n$
whenever we think of differential forms.
Twisting the Euler exact sequence  (cf.\,\cite[Thm 8.13, p.\,176]{Hart}):  we get
\begin{equation*}
0\to \Omega^1_{\p n}(d+2)\to \OX_{\p
  n}(d+1)\otimes S_{1}\to \OX_{\p n}(d+2)\to 0.
\end{equation*} 
 Taking global sections we find
 the exact sequence
\begin{equation}\label{Vdn}
0\to V_d^n:=H^0(\p n,\Omega^1_{\p n}(d+2))\lar
S_{d+1}\otimes S_{1}\stackrel{\iota_{R_{n}\vphantom{_|}}}\lar
S_{d+2}\to 0  
\end{equation} 
where $\iota_{R_{n}}(\sum A_{i}dZ_{i})=\sum A_{i}Z_{i}$ is the
contraction by the radial vector field.
Thus a 1-form $\w\in V_d^n$ can be 
written  in homogeneous coordinates as 
$$\w=A_{0}dZ_{0}+\dots +A_{n}dZ_{n}$$ where the $A_{i}$'s are
homogeneous polynomials of degree $d+1$ satisfying  
$$
A_{0}Z_{0}+\dots+A_{n}Z_{n} =0.
$$

\eco

\coi{Linear Pullback}\label{pullback}Let $F: \p n\dashrightarrow \p 2$ be a
  linear projection, \ie $F=[F_{0}:F_{1}:F_{2}]$,
  with $F_{i}$ linearly independent homogeneous
  polynomials of degree $1$. Pick $\w\in H^0(\p
  2,\Omega^1_{\p 2}(d+2))$. 
Write $\w=B_{0}dX_{0}+B_{1}dX_{1} +B_{2}dX_{2}$ where
  $B_{i}\in S_{d+1}(X_{0},X_{1},X_{2})$.  The pullback $F^*(\w)$ is  the 1-form
\begin{equation}\label{Fw}
F^*(\w)=F^*B_{0}dF_{0}+F^*B_{1}dF_{1} +F^*B_{2}dF_{2}.
\end{equation}
A simple application of Euler relation shows that
$\iota_{R_{n}}(F^*(\w))=F^*(\iota_{R_{2}}(w))=0$. On
the other hand, as any projective 1-form in $\p 2$
is integrable we have $0=F^*(\w\wedge
d\w)=F^*(\w)\wedge dF^*(\w)$ \ie
$F^*(\w)$ is a projective  integrable 1-form of degree $d$,
$$F^*(\w) \in H^0(\p n,\Omega^1_{\p n}((d+2))).$$   
\eco

\coi{Remark} 
For a fixed map $F$ as above we obtain injective
linear maps 
$$\ba c
H^0(\p 2,\Omega^1_{\p 2}(d+2))\lar
H^0(\p n,\Omega^1_{\p n}(d+2)),
\\\na7
\mathbb{F}(d,2) \lar \mathbb{F}(d,n).
\ea
$$
\eco

\coi{Parameter space for rational maps} Let  $F=[F_{0}:F_{1}:F_{2}]$ be a rational map as above,  with $F_{i}\in S_{1}$. 
Note  that if we change the basis of the linear system $\langle F_{0},F_{1},F_{2}\rangle$ then we obtain projectively equivalent pullbacks.  
So the  natural parameter space for rational
linear maps  is  the Grassmannian of dimension 3
subspaces of the space $S_1$ of forms of degree 1, 
\be\label{gr}
  \grass:=\grass(3,S_{1}).
\ee 

\eco

\section{Linear Pullback component}

We show next that the locus in $\p{}(V_d^n)$
(cf.\,\ref{Vdn}) corresponding to 
codimension one foliations obtained by linear
pullbacks of foliations in $\p 2$  is 
the birational image of a natural projective
bundle over the Grassmannian $\grass(3,S_1)$.

\begin{coisa}{\bf Proposition.}\label{existsfbgrass}
Notation as in {\rm(\ref{Vdn})}
and {\rm(\ref{gr})}
, let

\medskip
\centerline{$
\E_{d}:=\{(F,\mu)
\in  \grass\times V_d^n
\mid \mu=F\us\omega  \text{ \,for some  }
\,\omega\in V_d^2 
\} $.
}
\mn
 Then
\\ {\bf(i)} $\E_{d}$ is a vector subbundle
of \,$\grass\times V_d^n$\, of rank $(d+1)(d+3)$.
\\Let
$
q_2:\p{}(\E_{d})\subset\grass\times\p{}(V_d^n)
\to \p{}(V_d^n)$
  be the projection and set
  \\\centerline{$
\Y:=q_2(\p{}(\E_{d})), \ 
g:=\dim \grass=3(n-2).
$}
Then
\\{\bf(ii)} 
the dimension of  $\Y$ is \,
$g+(d+1)(d+3)-1$ \,
and
\\{\bf(iii)} 
the  degree of  $\Y$ is given by 
the Segre class  \  $
s_{g}(\E_{d}).
$ 
\end{coisa}

\begin{proof}

Consider  the tautological exact sequence of
vector bundles over 
$\grass$,  (cf.\,\ref{gr})

\begin{equation}\label{tau}
0\to \T\to \grass\times S_{1}\to \mathcal{Q}\to 0
\end{equation}
where $\T$ is of rank $3$ with fiber
$\T_{F}=\langle F_{0}, F_{1}, F_{2}\rangle, F_i\in
S_1$. 
We  obtain a natural  rational  map,
\begin{equation*}
\xymatrix
{\,\grass\times \p n\ar@{-->}[r]
^{ \overline F}
&  \p{}(\T^{\vee})
\ \ar@{->}[d]^{\psi} 
\\
&\grass
\ \ar@{<-}[ul] 
^{p_{1}} 
}
\end{equation*}
We 
interpret  the  map $\overline F$ as the family of
$F$'s as in \ref{pullback} and  $\p{}(\T^{\vee})$ as the family of $\p2$'s obtained as the fibers of  $\psi$. 
The relative  cotangent bundle of  $\psi$
fits into the relative Euler exact sequence (cf.\,\cite[B.5.8, p.\,435]{FUL}):  
\begin{equation*}\label{relcot}
0 \rightarrow\Omega_{\psi}(1)\rightarrow \psi^*\T\rightarrow \OX_{\T^{\vee}}(1)\rightarrow 0.
\end{equation*} 
Twisting by $ \OX_{\T^{\vee}}(d+1)$ we find
\begin{equation*}
0 \rightarrow\Omega_{\psi}(d+2)\rightarrow \psi^*\T(d+1)\rightarrow \OX_{\T^{\vee}}(d+2)\rightarrow 0.
\end{equation*} 
Taking direct image yields  the following exact
sequence over $\grass$ 
\begin{equation}\label{defVd}
0 \rightarrow \psi_{*}\Omega_{\psi}(d+2)\rightarrow \s{d+1}\T\otimes \T\rightarrow \s{d+2}(\T)\rightarrow 0.
\end{equation} 
Define 
\begin{equation} \label{Ednovo}
\E_{d}:=\psi_{*}\Omega_{\psi}(d+2).
\end{equation}
The fiber of $\E_{d}$  over each $F\in \grass$, 
 is   $$H^0(\p{}(\T_{F}^{\vee}),\Omega_{\p{}(\T_{F}^{\vee})}(d+2))$$ 
\ie  the
space of 
$1-$forms defining foliations of degree $d$  in
the varying  
$\p2\,\underline\sim\,\p{}(\T_{F}^{\vee})$.

On the other hand, we obtain from (\ref{tau})
{\em injective} maps of vector bundles over $\grass$,  
$$\iota_{1}:\s{d+1}\T\otimes \T\lar  S_{d+1}\otimes S_{1}$$ and 
$$\iota_{2}:\s{d+2}\T \lar  S_{d+2}.$$
These two maps fit into  the following  diagram of exact
sequences 
 
 \[
\xymatrix
{&\Ker (a) =\E_{d}\ar@{->}[r]^{j} \ar@{->}[d]&  \Ker (b) =V_{d}^n \ar@{->}[d]\\
0 \ar@{->}[r] & \s{d+1}\T\otimes \T \ar@{->}[r]^{\iota_{1}}\ar@{->}[d]^{a} & S_{d+1}\otimes S_{1}\ar@{->}[d]^{b}  \\
0 \ar@{->}[r]  &\s{d+2}\T\ar@{->}[r]^{\iota_{2}}& S_{d+2}\cdot
}
\]
In this way we obtain an injective map of vector bundles
over $\grass$,
\begin{equation}\label{jota}
j: \E_{d}\longrightarrow \grass \times V_{d}^n.
\end{equation}
The vector subbundle $
\E_{d}
\subset\grass\vez
V_{d}^n
$ is as stated in (i).

 Let $q_1:\p{}(\E_{d})\ra
\grass$ and $q_2:\p{}(\E_{d})\ra\p{}(V_{d}^n)$ be the 
maps induced by projection: 
\[
\xymatrix 
{
&\p{}(\E_{d}) 
 \ar[dl]_{q_{1
}} \ar[dr]^{q_{2}}&  \\
\grass & &
\ \ \ \ \ 
\Y\subset\p{}(V_{d}^n)}.
\] 
We prove in Lemma~\ref{GenInj} below  that  $q_2:\p{}(\E_{d})\ra\p{}(V_{d}^n)$ is
generically injective.
Set
$u:=\dim\Y=\dim\p{}(\E_{d})$. Write  $H$ for the
hyperplane class of $\p{}(V_{d}^n).$   We have
$q_2^*H=c_1\OX_{\E_{d}}(1)=:h$.
Using
\cite[\S3.1,\,p.\,47,\,Prop.4.4,\,p.\,83 and
Ex.\,8.3.14,\,p.\,143]{FUL}, we may compute
$$
\deg\Y=\int_{\p{}(V_{d}^n)}H^u\cap\Y=\int_{\p{}(\E_{d})}
h^u=\int_\grass q_{1*}(h^u)=\int_\grass{}
s_{g}(\E_{d}).
$$

\end{proof}

\begin{coisa}{\bf Remark} 
From sequence (\ref{defVd}) we obtain
$s(\E_{d})=s(\s{d+1}\T\otimes \T)c(\s{d+2}(\T))$.
\end{coisa}

\begin{coisa}{\bf Lemma.}\label{GenInj}
For general $(F,\w)\in \grass\times H^0(\p 2,\Omega^1_{\p 2}(d+2))$, $F^*(\w)$ determines uniquely  $\w$ and $F$.
\eco

\begin{proof}
With notation as in  \ref{pullback},  
let $I(F)=Z(F_{0},F_{1},F_{2})\subset \p n$ denote
the indeterminacy locus of $F$.   
For generic $(F,\w)$, the singular set of $F^*(\w)$ consists of linear
components of codimension two  of the form $F^{-1}(q)$
where $q\in Sing(\w)$. These linear components
intersect in  $I(F)$. Indeed, we can suppose $q_1=[0:0:1]$ and $q_2=[0:1:0]$, so $F^{-1}(q_1)=Z(F_{0},F_{1})$ and   $F^{-1}(q_2)=Z(F_{0},F_{2})$.
From $I(F)$ we retrieve the linear system $\langle
F_0,F_1,F_2\rangle$ \ie the point $F\in \grass$.

On the other hand, consider 
the blow-up $\pi: \mathbb{B}\lar \p n$
of $\p n$ in $I(F)$. We proceed to
show that it is possible to recover the 1-form  $\w$ from the strict transform $\tilde{\w}$ of $F^*(\w)$.
Indeed, recall
$\mathbb{B}=\{(p,[x_{0}:x_{1}:x_{2}])\in \p
n\times \p 2\mid
x_iF_j(p)-x_jF_i(p)=0\}$. Therefore   over $U:=\{x_{0}\neq 0\}$ we have 
$\mathbb{B}_{|U}=\{(p,[1:t:s])\mid F_1(p)=tF_0(p);
F_2(p)=sF_0(p)\}$ and the equation of the exceptional divisor is $F_0=0$.
In this chart we have $dF_1=F_0dt+tdF_0$,  $dF_2=F_0ds+sdF_0$, therefore 
\[\pi^*F^*(\w)=B_{0}(F_0,tF_0,sF_0)dF_{0}+B_{1}(F_0,tF_0,sF_0)dF_{1} +B_{2}(F_0,tF_0,sF_0)dF_{2}=\]
\[F_{0}^{d+1}[B_{0}(1,t,s)dF_{0}+B_{1}(1,t,s)(F_0dt+tdF_0) +B_{2}(1,t,s)(F_0ds+sdF_0)]=\]
\[F_{0}^{d+1} [(B_{0}(1,t,s)+tB_{1}(1,t,s) +sB_{2}(1,t,s))dF_0+B_{1}(1,t,s)F_0dt +B_{2}(1,t,s)F_0ds]\]

Recalling  $\w$ is a projective 1-form, so  $B_{0}(1,t,s)+tB_{1}(1,t,s) +sB_{2}(1,t,s)=0$,
we obtain 
\[\pi^*F^*(\w)=F_{0}^{d+2}[B_{1}(1,t,s)dt +B_{2}(1,t,s)ds].\]
Therefore  the strict transform of $F^*\w$ is \ 
$ \tilde{\w}=B_{1}(1,t,s)dt +B_{2}(1,t,s)ds.
$ \ 
Ditto for the other local charts of the blowup.
This shows that we may recover $\w$ from $\tilde{\w}$.

\end{proof}

\begin{coisa}{\bf Corollary.}
The component
$LPB(d,n)\subset\mathbb{F}(d,n)$ is rational. 
\eco

\section{Computation of the degree of  $LPB(d,n)$}


\begin{coisa}{\bf Proposition.}\label{degpolynomial}
Notation as in {\rm\ref{existsfbgrass}},   $\deg(LPB(d, n))$ is a polynomial in $d$ of degree $3g= 9(n-2)$.

 \end{coisa}
 
 We need some preliminary results.

\begin{coisa}{\bf Lemma.}\label{ck}
Let $E$ be a vector bundle of rank $r$ on a
variety $X$. The 
$k-$Segre class $$s_k(\s{d}(E))=\sum_{|\lambda|=k} p_\lambda(d) c_{\lambda}(E)$$ where $p_{\lambda}(d)$ is a polynomial  in $d$ of degree $\leq  r k$, and there exists $p_{\lambda}$ of degree $rk$.
Here the sum runs over the partitions of $k$,  and if $\lambda=(\lambda_1,\dots,\lambda_l)$ is a partition of $k$,  $c_{\lambda}(E):=c_{\lambda_1}(E)\cdots c_{\lambda_l}(E)$.

 \end{coisa}
\begin{proof}
 First we prove the following assertion relating
 the Chern characters to the Segre classes  of
 $\s{d}(E)$. 
\begin{coisa}{\bf Claim A(k):}\label{claim1}
 Assume the Chern character
$$ch_j(\s{d}(E^{\vee}))=\sum_{|\mu|=j} q_\mu(d) c_{\mu}(E)$$ where $q_{\mu}(d)$ are polynomials in $d$ of degree
 $\leq r+j-1$ for all $j\geq 0$.
 Then $s_k(\s{d}(E))=\sum_{|\lambda|=k} p_\lambda(d) c_{\lambda}(E)$ where $p_{\lambda}(d)$ are polynomials  in $d$ of degree $\leq  r k$. Moreover, if $ch_1(\s{d}(E))=q_1(d)c_1(E)$ and $deg(q_1)=r$, then $s_k(\s{d}(E))$ has a coefficient of degree $rk$.
\end{coisa}

  For a vector bundle $F$ of rank $r$, let
 $x_1,\dots,x_r$ be the Chern roots of
 $F^{\vee}$. Then for all $k\geq 0$, the Segre
 class \[s_k(F)=\sum_{1\leq i_1\leq \dots\leq
 i_k\leq r}x_{i_1}\cdots x_{i_k}\] is the
 $k-$\emph{complete symmetric function}
 (cf. \cite[p.\,28]{FP}).


On the other hand the complete symmetric functions
can be expressed in term of the \emph{power sum
symmetric functions}, $p_k=\sum_{i} x_i^k$. We borrow from
(\cite[ p.\,25]{Mac}) the explicit relations:
\be
s_k(F)=\sum_{|\lambda|=k} w_{\lambda} ch_{\lambda}(F^{\vee})
\label{CHCH}
\ee
where $\lambda=(\lambda_1,\dots,\lambda_l )$ is a partition of $k$. Following the notation in \cite{Mac} we write $\lambda=(1^{m_1},2^{m_2}, \dots)$ where $m_i:=\#\{j\mid \lambda_j=i\}$ and 
$w_\lambda=\prod_{i\geq 1}\frac{\lambda_i!}{i^{m_i}m_i!}$.

Write 
 for short $ch_j=ch_j(\s{d}(E^{\vee}))$.
Whenever the coefficients of
 $ch_j$ are  polynomials in $d$ of
 degree $m$ we will write
 $\deg(ch_j) = m$,  by abuse of notation.  
 We observe that $deg(ch_{\lambda})=deg(ch_{\lambda_1}\cdots
ch_{\lambda_l})\leq l(r-1)+k\leq rk$ and the equality holds if and only if $l=k$, i.e. if $\lambda=(1,\dots 1)$ in which case the coefficient of $ch_1^k$ is $\frac{1}{k!}$. In other words \[s_k(\s{d}(E))=\frac{1}{k!}ch_1^k+\text{l.o.t}\]
Hence $s_k(\s{d}(E))$ is a linear
 combination of  monomials in the  Chern classes
 of $E$ whose    
 coefficients are polynomials in $d$ of degree
 $\leq rk$.

Next we prove the following claim by induction on $r=\rk(E)$ and on  $k$:

\vspace{.2cm}
\begin{coisa}{\bf Claim P(k):}
\label{claim}
For $k\geq 0$, $ch_k(\s{d}(E))=\sum_{|\mu|=j} q_\mu(d) c_{\mu}(E)$, where $q_{\mu}(d)$ are polynomials in $d$ of degree
 $\leq r+k-1$.
  Moreover, $ch_{1}(\s{d}(E))=q_1(d)c_1(E)$ where $q_1(d)$ is a polynomial of degree $r$.
\end{coisa}
\vspace{.2cm}
For $r=1$ we have $ch_{k}(\s{d}(E))=\frac{1}{k!}d^kc_1(E)^k$. 

Suppose that \mb{P(k)} is true for vector bundles of rank $r-1$. Let $\pi: E\to X$ be a vector bundle of $\rk(E)=r$ and $p:\p{}(E) \to X$ the induced projective bundle.
For $k=0$ we have $ch_0(\s{d}(E))=\rk(\s{d}(E))=\binom{d+r-1}{r-1}$, a polynomial in $d$ of degree $r-1$.
Suppose that  $k\geq 1$, and that \mb{P(s)} holds
for $\mb s<k$. 

Over $\p{}(E)$ we have the tautological exact sequence: 
\begin{equation}\label{tauE}
0\to \OX_{E}(-1)\to p^*E\to Q\to 0.
\end{equation} 
It  induces the following exact sequence 
 for $d\geq 1$
\begin{equation*}
0\to \OX_{E}(-1)\otimes \s{d-1}(p^*E)\to \s{d}(p^*E)\to \s{d}(Q)\to 0.
\end{equation*} 
Hence we may write the
relation for the Chern characters  
\begin{equation}\label{sumch}
ch(\s{d}(p^*E))=ch(\OX_{E}(-1)\otimes \s{d-1}(p^*E))+ch(\s{d}(Q)).
\end{equation}
On the other hand,
$ch\left(\OX_{E}(-1)\otimes \s{d-1}(p^*E)\right)=ch(\OX_{E}(-1)) 
ch(\s{d-1}(p^*E)$. So each graded part satisfies
\begin{equation}\label{OE}
ch_k\left(\OX_{E}(-1)\otimes \s{d-1}(p^*E)\right)=
\sum_{i=0}^{k}\frac{\alpha^i}{i!}ch_{k-i}\left
(\s{d-1}(p^*E)\right) 
\end{equation}
where $\alpha:=c_1(\OX_{E}(-1))$.
It follows from (\ref{sumch}) and (\ref{OE}) that
\begin{equation*}
ch_k(\s{d}(p^*E))=\sum_{i=0}^{k}\frac{\alpha^i}{i!}ch_{k-i}(\s{d-1}(p^*E))+ch_{k}(\s{d}(Q)).
\end{equation*}
Hence 
\begin{equation}\label{diff}
\ba r
ch_k(\s{d}(p^*E))-ch_k(\s{d-1}(p^*E))=\sum_{i=1}^{k}\frac{\alpha^i}{i!}ch_{k-i}(\s{d-1}(p^*E))
\\
+ch_{k}(\s{d}(Q)).
\ea
\end{equation}

Observe that the right hand side of (\ref{diff})
involves: 
\begin{itemize}
\item $ch_k(\s{d}(Q))$, which by  induction,
since $\rk(Q)=r-1$,  is of the form  $\sum_{|\mu|=k} q_\mu(d) c_{\mu}(Q)$ where $q_{\mu}(d)$ are polynomials
 in $d$ of degree  $\leq r-1 +k-1=r+k-2$.
 Moreover, by (\ref{tauE})
  $c_s(Q)=c_s(p^*E)-\alpha c_{s-1}(p^*E)$. Thus
  $ch_k(\s{d}(Q))$ is a linear  combination of monomials
  in the Chern classes of $p^*E$ and in $\alpha$  whose   
 coefficients are polynomials in $d$ of degree
 $\leq r+k-2$. 
 
\item  $ch_{s}(\s{d-1}(p^*E))$ with $s<k$ which,
by induction  on $k$ is of the form $\sum_{|\mu|=k} q_\mu(d) c_{\mu}(p^*(E))$ where $q_\mu(d)$ are polynomials in $d$ of degree
 $\leq r+s-1$;  the maximal  degree appearing is
 $ r+k-2$ (coming from the coefficients in  $\alpha
 ch_{k-1}(\s{d-1}(p^*E))$).   
\end{itemize}

Recall that the pullback $p^*: A_*(X)\to
A_*(\p{}(E))$ is a monomorphism  with left inverse  
$\beta\mapsto p_*(c_1(\OX_{E}(1))^{r-1}\cap (\beta))$ ( \cite[p.\,49]{FUL}).
Applying this inverse  to (\ref{diff}) we conclude
that
$$
ch_k(\s{d}(E))-ch_k(\s{d-1}(E))
$$ is a linear combination of monomials in the Chern classes of $E$ whose coefficients are polynomials in $d$ of degree $\leq r+k-2$, and this implies that the coefficients in $ch_k(\s{d}(E))$   are   polynomials in $d$ of degree $\leq r+k-1$. 

Observe that for $k=1$ we obtain 
that $ch_1(\s{d}(E))-ch_1(\s{d-1}(E))=\binom{d+r-2}{r-1}c_1(E)$, a polynomial of degree $r-1$. So $ch_1(\s{d}(E))=q_1(d)c_1(E)$ where $q_1(d)$ is polynomial of degree $r$.

Using \mb{A(k)} \ref{claim1}, we deduce  that $s_k(\s{d}(E))$
is a linear combination of monomials in the Chern classes of $E$  whose coefficients are polynomials in $d$ of degree $\leq \rk(E)k$ and there exists a coefficient of degree $\rk(E)k$.

\end{proof}

Next we prove Proposition \ref{degpolynomial}.

\begin{proof}

From the exact sequence
\begin{equation*}
0 \rightarrow \mathcal{V}_d \rightarrow \s{d+1}\T\otimes \T\rightarrow \s{d+2}(\T)\rightarrow 0
\end{equation*} 
  we obtain the relation for Chern characters
  $$ch(\mathcal{V}_d^{\vee})=ch(\s{d+1}\T^{\vee}\otimes \T^{\vee})-ch(\s{d+2}(\T^{\vee})).$$
 Thus for each $k\geq 1$ we have
 $$ ch_k(\mathcal{V}_d^{\vee})=\sum_{i=0}^k
 ch_{k-i}(\s{d+1}\T^{\vee})ch_i(\T^{\vee})-ch_k(\s{d+2}(\T^{\vee})).
 $$
  By  assertion  \mb{P(k)} (\ref{claim})
$ch_k(\s{d}(\T^{\vee}))$ is a polynomial in $d$ of degree $\leq k+2$. Therefore  $ch_k(\mathcal{V}_d^{\vee})$ is a polynomial in $d$ of degree $\leq k+2$. 
Moreover,
\[ch_{1}(\s{d}(\T^{\vee}))=q_1(d)c_1(\T^{\vee})
\]
where $q_1(d)$is a polynomial of degree $3$. Therefore
\[ch_1(\mathcal{V}_d^{\vee})=
 ch_{1}(\s{d+1}\T^{\vee})ch_0(\T^{\vee})+ch_{0}(\s{d+1}\T^{\vee})ch_1(\T^{\vee})-ch_1(\s{d+2}(\T^{\vee}))=\]
\[(3q_1(d+1)+\binom{d+3}{2}-q_1(d+2)) c_1(\T^{\vee})\] is a polynomial in $d$ of degree $3$.

As in the proof of  assertion $\mb{A(k)}$
(\ref{claim1}), we conclude that $s_k(\mathcal{V}_d)$ is polynomial in $d$ of degree $3k$.
 
\end{proof}

\coi{Some formulas}\label{formulas}
To get explicit formulas for
$s_{g}(\E_{d})$ for any fixed $d,n$ we use sequence 
(\ref{defVd}),
and Macaulay2,\,\cite{Macaulay2}. We find (cf.\,
script
\S\ref
{scripts} below)

$\deg(LPB(d, 3))=
\frac{20}{27}
\binom{d+4}5
(d^2 + 6d + 11)(d^2 + 2d + 3)$.

$\deg(LPB(d, 4)) =
\frac{1}{839808}\frac{(d+4)!}{(d-1)!} (8d^{12} + 192d^{11} + 2176d^{10} + 15360d^9+
75090d^8 + 267552d^7 + 711859d^6 + 1423716d^5+
2119892d^4 + 2279136d^3 + 1662291d^2 + 730188d+
125388)(2 + d)$.

\eco

 \section{Higher degree pullback components}

As stated in the Introduction, the set of
foliations obtained by pullback of foliations in
$\p 2$ by rational maps of degree $m$ also form an
irreducible component $PB(m, k,n)$ of
$\mathbb{F}(d,n),\,d:=(k+2)m-2$.

As in the linear  case, a natural parameter
space of rational maps
$F: \p n\dashrightarrow \p 2$   of degree $m$ is
the grassmannian
$
\grass(3,S_m)$. Mimicking \ref{existsfbgrass},
we can construct a fiber bundle

\centerline{$
\E_{m,k}:=\{(F,\mu)
\in  \grass(3,S_m)\times V_k^n
\mid \mu=F\us\omega  \text{ \,for some  }
\,\omega\in V_k^2 
\} $.
}
\noindent
However, for $m\geq2$, 
the  map  $j: \E_{m,k}\to V_{d}^n$ (analogous to the
map in (\ref{jota})) is injective only over the
open subset consisting of dominant maps. So its
image is not  a subundle of $V_{d}^n$. In fact,
for a non-dominant map $F\in \grass(3,S_m)$,
$j_F(\w)=F^*(\w)=0$ for all $\w\in H^0(\p
2,\Omega^1_{\p 2}(k+2))$ defining a foliation that
leaves invariant
the closure of
${\im(F)}$. 

We could in principle find the locus $Z\subset \grass(3,S_m)$ where the rank of $\im(j)$ drops and then  blow-up $\grass(3,S_m)$ along $Z$. Doing this we expect to build a subundle $\widetilde{\E}_k$ of $V_{d}^n$ over the blow-up that coincides with $\im(j)$ over the open set of dominant maps.  The  projection of $\p{}(\widetilde{\E}_k)$ to  $\p{}(V_d^n)$ is  the   space of foliations obtained as pullback by dominant maps.
To make it work, we'd need to know how to get our
hands on the Segre classes of $Z$.
We hope to report on this elsewhere.

\section{scripts}
\label{scripts}

\coi{Scripts for Macaulay2} In order to compute
$\deg(LPB(d, n))$
for any given value of $n\geq3$, 
just set $N=n$ at the beginning of
the script below. It can be fed into
\url{http://habanero.math.cornell.edu:3690/}.

{\small
\begin{verbatim}

loadPackage "Schubert2"
N=3  --plug-in 3,4...
pt= base d
-- set d to  be a free parameter in the ``intersection ring'' 
--of the base variety
G=flagBundle({3,N-2}, pt)
-- Grassmannian of 3-planes in N+1-space,
(S,Q)=G.Bundles
-- names the sub and quotient bundles on G
A=symmetricPower(d+2, S)
B=symmetricPower(d+1, S)*S
integral(chern(A-B))


\end{verbatim}
}

\eco

{\bf Acknowledgments}.
We wish to thank Luca Scala for  the  suggestions that  led to a significant improvement of the proof of Proposition \ref{degpolynomial}.

\end{document}